\documentclass[12pt, reqno]{amsart}
\usepackage{amsmath, amsthm, amscd, amsfonts, amssymb, graphicx, color}
\usepackage{setspace}
\usepackage{mathrsfs}
\usepackage{multicol}
\usepackage[bookmarksnumbered, colorlinks, plainpages]{hyperref}
\hypersetup{colorlinks=true,linkcolor=red, anchorcolor=green, citecolor=cyan, urlcolor=red, filecolor=magenta, pdftoolbar=true}

\newtheorem{theorem}{Theorem}[section]
\newtheorem{lemma}[theorem]{Lemma}
\newtheorem{proposition}[theorem]{Proposition}
\newtheorem{corollary}[theorem]{Corollary}
\theoremstyle{definition}

\theoremstyle{remark}

\numberwithin{equation}{section}

\newcommand{\NN}{\mathbb{N}}

\newcommand{\CC}{\mathbb {C}}

\newcommand{\RR}{\mathbb{R}}

\begin{document}
\setcounter{page}{1}
\title[ Surjective and  closed range differentiation  operator  ]{ Surjective and closed range differentiation  operator  }
\author[Tesfa  Mengestie]{Tesfa  Mengestie}
\address{ Mathematics Section \\
Western Norway University of Applied Sciences\\
Klingenbergvegen 8, N-5414 Stord, Norway}
\email{Tesfa.Mengestie@hvl.no}
\makeatletter
\@namedef{subjclassname@2020}{%
  \textup{2020} Mathematics Subject Classification}
\subjclass[2020]{Primary 47B37, 30H20,  Secondary 46E15, 46E20 \newline\hspace*{2em}\emph{Keywords and phrases.}  Fock-type  spaces; Closed range; Differentiation operator; Surjective; Order bounded}
 \begin{abstract}
  We  identify   Fock-type spaces $\mathcal{F}_{(m,p)}$  on which the differentiation  operator $D$ has   closed range. We prove that $D$ has closed range only if it is surjective, and this happens if and only if
  $m=1$.  Moreover, since the operator is unbounded on the classical Fock spaces, we consider the modified or the weighted   composition--differentiation operator, $D_{(u,\psi,n)} f= u\cdot\big( f^{(n)}\circ \psi\big)$,  on these spaces and  describe   conditions under which the operator  admits closed range, surjective, and order bounded structures.
 \end{abstract}
\maketitle
\section{Introduction}
The  differentiation  operator,   $Df= f'$,   often appears  as  an   example of unbounded linear operators  in many Banach spaces,  including Hardy spaces and  Bergman spaces  \cite{HP},   Fock spaces  and   Fock-type   spaces where the weight decays faster than the Gaussian weight  \cite{TM3}, and  Fock--Sobolev spaces,  where the weights decay slower than  the Gaussian weight \cite{TM5}. Inspired by all these developments,  the  question  whether there  could exist Fock-type spaces   on which  $D$ admits  basic operator-theoretic  structures  was  investigated in  \cite{TMD}.
To answer the question, the author considered  the space
\begin{align*}
\mathcal{F}_{(m,p)}:= \Big\{ f\in \mathcal{H}(\CC) :\|f\|_{(m, p)}^p= \int_{\CC} |f(z)|^p e^{-p|z|^m} dA(z) <\infty\Big\},
\end{align*} where   $\mathcal{H}(\CC)$ denotes the set of entire functions on the complex plane $\CC$,   $m>0$,   $1\leq p<\infty$,  and $A$ is  the
usual Lebesgue area  measure.   Then the following basic property  was proved.
\begin{theorem}\label{thmDD}\cite[Theorem 1.1]{TMD} Let $1\leq p, q< \infty$ and $m>0$.
\begin{enumerate}
\item  If $p\leq q$, then  $D: \mathcal{F}_ {(m,p)} \to \mathcal{F}_{(m,q)}$ is
 bounded if and only if \begin{align}
 \label{diffin}
 m\leq   2-\frac{pq}{pq+q-p},\end{align}
 and compact if and only if  the inequality in \eqref{diffin}  is strict.
\item If $p>q$, then  $D: \mathcal{F}_{(m,p)} \to \mathcal{F}_{(m,q)}$ is bounded (compact) if and only if
 \begin{align*}
 m<  1-2  \Big(\frac{1}{q}-\frac{1}{p}\Big).
\end{align*}
\end{enumerate}
\end{theorem}
For $p= q$, the inequality in \eqref{diffin} simplifies  to  $m\leq 1$,  which is stronger than the boundedness  condition for $p<q$. On the  corresponding growth type space
\begin{align*}
\mathcal{F}_{(m, \infty)}:=\big\{ f\in \mathcal{H}(\CC): \|f\|_{(m, \infty)}= \sup_{z\in \CC} |f(z)|e^{-|z|^m}<\infty \big\},
\end{align*}
the boundedness of $D$ was characterized by the same  condition $m\leq 1$; see  \cite{Maria,Bon, Harutyunyan}.
If $D$ acts between two  different Fock-type spaces $\mathcal{F}_{(m,p)}$ and $\mathcal{F}_{(m,q)}$,  where one of the spaces is growth type, then a simple variant of the proof of Theorem~\ref{thmDD} in \cite{TMD} gives the following result.
\begin{corollary}\label{corr}
Let $1\leq p< \infty$ and $m>0$. Then  the operator
\begin{enumerate}
\item    $D: \mathcal{F}_ {(m,p)} \to \mathcal{F}_{(m,\infty)}$ is
 bounded if and only if \begin{align}
 \label{diffin1}
 m\leq   2-\frac{p}{p+1},\end{align}
 and compact if and only if  the inequality in \eqref{diffin1}  is strict.
\item  $D: \mathcal{F}_{(m,\infty)} \to \mathcal{F}_{(m,p)}$ is bounded (compact) if and only if
$ m< 1-\frac{2}{p}.$
\end{enumerate}
\end{corollary}
For more related results, we   refer the   interested readers to \cite{BMM,TMD} and the references therein. One of the main objectives of this work is to  identify Fock-type spaces on which the  differentiation operator  admits  closed range structure.   Our next  main result   shows   there exists no closed range differentiation  operator acting between two different Fock-type spaces.
\begin{theorem} \label{closeddiff1}
Let $1\leq p, q \leq \infty$, $m>0$,  and   $D:\mathcal{F}_{(m,p)} \to  \mathcal{F}_{(m,q)}$ be bounded.  Then the following statements are equivalent.
\begin{enumerate}
\item   $D$  has closed range;
\item $p=q $ and $m=  1$;
\item $D$   is  surjective.
\end{enumerate}
\end{theorem}
The result identifies $\mathcal{F}_{(1,p)}$ as the only Fock-type space  supporting  closed range structure for the operator $D$.
The proof of the result will be presented in Section~\ref{mainproof}.

  As mentioned earlier, $D$ is not  bounded  on the classical Fock spaces. In \cite{TMDV}, the   author studied   whether simply modulating the classical Gaussian weight function $|z|^2/2$ by positive parameters  $\alpha$ would produce a bounded  $D$ on the  spaces
\begin{align*}
\mathcal{F}_\alpha^p:=\Big\{ f\in \mathcal{H}(\CC): \int_{\CC} |f(z)|^p
e^{-\frac{p\alpha}{2}|z|^2}  dA(z)< \infty\Big\}.
\end{align*} It follows that for positive parameters $\alpha$ and $\beta$,   $D: \mathcal{F}_{\alpha}^p \to \mathcal{F}_{\beta}^q$ is bounded if and only if $\alpha  <\beta $  and $p\leq q$. This condition equivalently describes the compactness of
the operator.  Consequently, as it will  be  explained     later in  the proof of Theorem~\ref{closeddiff1}, such  modulated  spaces support no closed range compact differentiation  operator.
\subsection{Weighted composition--differentiation  operator}
In the preceding section, we considered Fock-type spaces on which  the differentiation  operator is bounded,  and  we classified them  based on whether they support
 closed range structure for $D$ or not. In this section,  we modify  the operator itself and study the closed range and surjectivity  problems on the classical Fock spaces. \\
  For each $n\in \NN_0=\{0, 1, 2, ...\}$ and   $u , \psi$ in   $\mathcal{H}(\CC)$, we  define
the weighted  composition--differentiation operator  $ D_{(u,\psi,n)} $ by
 \begin{align*}
D_{(u,\psi,n)} f= u\cdot\big( f^{(n)}\circ \psi\big),
\end{align*} where $f^{(n)}$ is the $n^{th}$ order derivative of the function $f$ and $f^{(0)}= f$. Then we investigate  when the operator $D_{(u,\psi,n)} $ admits closed range structure  on the classical  Fock spaces.
   This class of operators has lately  attracted a considerable amount of attention; see for example  \cite{MW} and the references therein.\\
  For  $1\leq p\leq \infty$, recall that the   Fock spaces $\mathcal{F}_p$  are defined by
  \begin{align*}
\mathcal{F}_p:=\big\{f\in \mathcal{H}(\CC): \|f\|_{p} <\infty \big\},
\end{align*} where
\begin{equation*}
\|f\|_{p}:=
\begin{cases}\bigg(\frac{ p}{2\pi}\int_{\CC} |f(z)|^p
e^{-\frac{ p}{2}|z|^2} dA(z)\bigg)^{\frac{1}{p}} <\infty, \ \ 1\leq p<\infty \\
\sup_{z\in \CC}
|f(z)|e^{-\frac{1}{2}|z|^2} <\infty, \ \ p= \infty.
\end{cases}
\end{equation*}  The space  $\mathcal{F}_2$ is a
reproducing kernel Hilbert space with kernel and normalized reproducing kernel
functions  \begin{align*} K_{w}(z)= e^{ \overline{w}z}\ \   \text{and} \ \
k_{w}(z)=\|K_w\|_2^{-1}K_w(z)= e^{\overline{w}z-\frac{|w|^2}{2}}
\end{align*} for all $z, w\in \CC$.  A straightforward  calculation shows   $\|k_w\|_p= 1 $ for all $p$. For more details, we refer to the book \cite{Zhu}.

 Note that for each   $f\in  \mathcal{H}(\CC)$ and  $p\neq \infty$, by \cite[p. 37]{Zhu} the local estimate
 \begin{align}
 \label{localest}
 |f(z)|
 \leq \frac{e^{\frac{|z|^2}{2}}}{r^2}\bigg( \int_{D(z,r )} |f(w)|^p
e^{-\frac{p|w|^2}{2}} dA(w)\bigg)^{1/p}
 \end{align} holds, where $D(z, r)$ is a disc of center $z$ and radius $r$. For $r=1$, this  estimate gives
   \begin{align}
 \label{localest2}
 |f(z)| \leq  e^{\frac{|z|^2}{2}}\|f\|_p.
 \end{align}
 By definition of the norm, the estimate in \eqref{localest2} holds for $p= \infty$  as well.

We note that the result in \cite{TM3} addresses  the unboundedness of only  the  first order differentiation operator $D$ on  Fock spaces.  A simple argument shows the $n^{th}$ order differentiation  operator, $D^n f= f^{(n)}$, is not bounded  for all $n\in \NN$ either. Indeed, using the kernel function $K_w$, we observe
\begin{align*}
\frac{\|D^n K_w\|_q}{\|K_w\|_p}= |\overline{w}|^n \frac{\|K_w\|_q}{\|K_w\|_p}= |w|^n \to \infty,
\end{align*} when  $|w|\to \infty$ independently of the exponents $p$ and $q$. On the other hand, an easy computation
using the equivalent norms in \eqref{HU} below shows  the composed differentiation  operator, $D^n C_\psi f= f^{(n)}\circ \psi$, is bounded on Fock spaces for every $\psi(z)= az$  and  $ 0<|a|<1$. Motivated by this, in \cite{MW} the bounded and compact structures of   $D_{(u,\psi, n)}$  were described in terms of the function
    \begin{align*}
 L_{(u,\psi,n)}(z):= |u(z)||\psi(z)|^n e^{\frac{1}{2}( |\psi(z)|^2- |z|^2)}
  \end{align*}for each $z \in \CC$. For further referencing, we state the result below.
    \begin{theorem}\cite[Theorem 1.3]{MW} \label{thm3}
  Let   $u, \psi\in \mathcal{H}(\CC)$, $n\in \NN_0$,  and $1\leq p, q\leq  \infty$.
  \begin{enumerate}
  \item If $p\leq q$,  then  $D_{(u,\psi,n)}: \mathcal{F}_p \to \mathcal{F}_q$ is bounded  if and only if   $L_n:=  \sup_{z\in \CC} L_{(u,\psi, n)}(z) <\infty$, and
   compact  if and only if $L_{(u,\psi,n)}(z) \to 0$ when $|z| \to \infty$.
   \item If $p> q$,  then  $D_{(u,\psi,n)}: \mathcal{F}_p \to \mathcal{F}_q$ is bounded (compact) if and only if $ L_{(u,\psi,n)} \in L^\frac{pq}{p-q}(\CC, dA)$ for $p<\infty$ and
   $ L_{(u,\psi,n)} \in L^q(\CC, dA)$ for $p= \infty$.
  \end{enumerate}
  \end{theorem}
    A bounded   $D_{(u,\psi, n)}$ implies
    $\psi(z)= a_nz +b_n$,  $a_n, b_n \in \CC$ and $|a_n|\leq 1$: see \cite {MW} for the details. For simplicity, we   write  $\psi(z)= az+b, a, b \in \CC$.   If $|a|=1$,  then
    \begin{align*}
 L_{(u,\psi,n )}(z)= |u(z)||\psi(z)|^n e^{\frac{1}{2}( |\psi(z)|^2- |z|^2)}
   =  |u(z)\psi^n(z) K_{\overline{a}b}(z)| e^{\frac{|b|^2}{2}}  \leq L_n
    \end{align*}  for all $z\in \CC$, $n\in \NN_0$,  and $L_n$ as in Theorem~\ref{thm3}.  Consequently,
    \begin{align*}
    |u(z)\psi^n(z) K_{\overline{a}b}(z)|   \leq L_n e^{-\frac{|b|^2}{2}}.
    \end{align*}
    By Liouville's theorem, it follows  that   $u \psi^n K_{\overline{a}b}$ is a constant $C_n$ and hence
    $u(z) \psi^n(z)= C_n K_{-\overline{a}b}(z)$. Setting $ z=0$, we get  $ C_n= b^nu(0)$. Therefore,
    \begin{align}
    \label{kernform}
    u(z) \psi^n(z)= b^nu(0) K_{-\overline{a}b}(z).
    \end{align} The representation in \eqref{kernform} will   be  needed later.

  For an  $f\in \mathcal{H}(\CC)$ and a positive $r$, we set $M_f(r)=\max\{ |f(z)|: |z|= r\}$. Then the order $\rho(f)$ of $f$ is defined by
\begin{align*}
\rho(f)= \limsup_{r\to \infty} \frac{\log\big(\log M_f(r)\big)}{\log r}.
\end{align*}
Now by Theorem~\ref{thm3},
\begin{align}
|u(z)||\psi(z)|^n  \leq L_n e^{\frac{1}{2}(|z|^2-|az+b|^2)}
\end{align} for all $z\in \CC$. It follows that $ u\psi ^n$ is of  order at most $ 2$.  Consequently,  a simple variant of the proof of  \cite[Theorem 3.2]{CG} gives
the  following representation  whenever $u\psi^n$ is non-vanishing.
 \begin{lemma} \label{lemnew}
 Let $u, \psi\in \mathcal{H}(\CC)$, $n\in \NN_0$, and  $1\leq p\leq  \infty$. Let  $D_{(u,\psi, n)}$ be bounded on $\mathcal{F}_p $  and hence $\psi(z)= az+b$  for some $a, b \in \CC$  such that $|a|\leq1$.
If $0\le |a|< 1$ and  $u\psi^n$ is non-vanishing, then
    \begin{enumerate}
 \item  $D_{(u,\psi, n)}$ is  compact  if and only if $ u\psi^n$ has the form
  \begin{align}
  \label{neww}
  u(z)\psi^n(z)= e^{a_{0n}+a_{1n}z+ a_{2n}z^2}
    \end{align} for some constants $a_{0n}, a_{1n}, a_{2n} $ in $ \CC$ such that $|a_{2n}|<\frac{1-|a|^2}{2}$.
    \item  $D_{(u,\psi, n)}$ is   bounded but not compact if and only if
    $u\psi^n$ has the form in \eqref{neww} with  $|a_{2n}|=\frac{1-|a|^2}{2}$ and either   $a_{1n}+a\overline{b}=0$ or $a_{1n}+a\overline{b}\neq 0$ and $$a_{2n}= - \frac{(1-|a|^2)(a_{1n}+a\overline{b})^2}{2|a_{1n}+a\overline{b}|^2}. $$
   \end{enumerate}
\end{lemma}
 By Theorem~\ref{thm3} we get that for $u(z)=1$ and $\psi(z)= az+b$  with  $|a|<1$   the operator $D_{(1, \psi, n)}$  is  compact.   In this case, the operator is
quasinilpotent,  that is, its spectral radius is zero.   Indeed, if $\lambda$ is an eigenvalue, then
\begin{align*}
D_{(1, \psi, n)} f=  f^{(n)}(az+b)= \lambda f(z)
\end{align*} for some  nonzero function
$f\in\mathcal{F}_p $.
Differentiating $m$ times both sides of the equation gives
\begin{align*}
 a^m f^{(n+m)}(az+b)= \lambda f^{(m)}(z).
\end{align*}
It follows that
\begin{align*}
  f^{(n+m)}(az+b)=  D_{(1, \psi, n)}) f^{(m)}(z)= \frac{\lambda}{a^m} f^{(m)}(z).
\end{align*}
This shows  $f$ cannot be a polynomial and $f^{(m)}$ cannot be zero. Thus, $\lambda/a^m$ forms a sequence of eigenvalues for  $D_{(1, \psi, n)}$,  which is a contradiction since
 $\lambda/a^m \to \infty $ as $m\to \infty$.
\subsection{Order bounded  $D_{(u,\psi, n)}$}
Another  interesting notion  closely related to  boundedness of an  operator is  order boundedness. The notion finds applications due to its close  relation to   absolutely  summing operators \cite{JDHJ,TD}.
 We say that an operator  $T:\mathcal{F}_p\to \mathcal{F}_q$ is order bounded  if there exists a positive function $h\in L^q\big(\CC, e^{-\frac{q}{2}|z|^2}dA(z)\big) $  such that for all $f\in \mathcal{F}_p$ with  $\|f\|_p \leq 1$
\begin{align*}|T(f(z))| \leq h(z)\end{align*}
almost everywhere with  respect to the measure $A$.  This  definition  was introduced by Hunziker and Jarchow in \cite{HHJ}. For some recent work on the subject, we refer the interested reader
to \cite{AKS,EW} and the references therein.  
  For  the operator $D_{(u,\psi,n)}$, we provide the following  characterization.
\begin{theorem}
\label{thmorder}
Let   $u, \psi\in \mathcal{H}(\CC)$, $n\in \NN_0$,  and $1\leq p, q\leq  \infty$. Then
 $D_{(u,\psi,n)}: \mathcal{F}_p \to \mathcal{F}_q$ is order bounded   if and only if   $ L_{(u,\psi, n)}\in L^q(\CC, dA)$.
   \end{theorem}
     By the discussion after  Theorem~\ref{thm3}, we observe  that the order bounded condition implies  $\psi(z)= az+b$ and $|a|<1$.  Note that if $|a|=1,$ then by \eqref{kernform}
   \begin{align*}
 L_{(u,\psi,n )}(z)=  |u(z)\psi^n(z) K_{\overline{a}b}(z)| e^{\frac{|b|^2}{2}} =| b^nu(0)| |K_{-\overline{a}b}(z) K_{\overline{a}b}(z)| e^{\frac{|b|^2}{2}}=| b^nu(0)|e^{\frac{|b|^2}{2}}
    \end{align*} and hence the condition in the theorem fails to hold.  On the other hand, for $ q= \infty$, the boundedness and order boundedness conditions
    coincide.  This together with Theorem~\ref{thm3} implies   every order bounded operator $D_{(u,\psi,n)}$ is compact.  By \cite[Theorem~1.8]{MW},  we also observe that the order bounded $D_{(u,\psi,n)}$  are  exactly
    those which are Hilbert-Schmidt in  $\mathcal{F}_2$.

     If $n=0$, then  $D_{(u,\psi, n)}$ reduces to the weighted composition operator $D_{(u,\psi, 0)} f= u \cdot( f\circ \psi)=W_{(u,\psi)}$.  Thus, Theorem~\ref{thmorder} describes
   the order bounded weighted composition operators on Fock spaces, which we state it as follows.
    \begin{corollary}
\label{corr03}
Let   $u, \psi\in \mathcal{H}(\CC)$  and $1\leq p, q\leq  \infty$. Then
 $ W_{(u,\psi)}: \mathcal{F}_p \to \mathcal{F}_q$ is order bounded   if and only if   $ m_{(u,\psi)}\in L^q(\CC, dA)$,
  where
  \begin{align*}
  m_{(u,\psi)}(z)=|u(z)| e^{\frac{1}{2}( |\psi(z)|^2- |z|^2)}
  \end{align*} for all $z\in \CC$.
   \end{corollary}
  If  $u=1$, then  we set   $D_{(\psi,n)}:= D_{(1,\psi,n)} $, and  Theorem~\ref{thmorder} implies  the following result  about
  the composition-differentiation operator.
  \begin{corollary}
  Let   $ \psi\in \mathcal{H}(\CC)$  and $1\leq p, q\leq  \infty$. Then for each $ n \in \NN_0$, the following statements are equivalent.
  \begin{enumerate}
  \item $D_{(\psi,n)}: \mathcal{F}_p \to \mathcal{F}_q$ is bounded;
 \item  $D_{(\psi,n)}: \mathcal{F}_p \to \mathcal{F}_q$ is order bounded;
 \item   $\psi(z)= az+b$ and  $|a|<1$.
 \end{enumerate}
  \end{corollary}
 \subsection{ Closed range $D_{(u,\psi, n)}$}
 We now  study the closed range property of  $ D_{(u,\psi,n)}$.
If $\psi= b \in \CC$, then  $D_{(u,\psi,n)}f= u f^{(n)}(b) $, and hence the range,
\begin{align*}
\mathcal{R}(D_{(u,\psi,n)})= \big\{u f^{(n)}(b):  f\in \mathcal{F}_p\big\},
\end{align*} is closed.  Thus, we assume  that $\psi$ is not a  constant in the rest of the manuscript. The next proposition  shows a nontrivial $D_{(u,\psi,n )}$ cannot have closed range if it acts between two different Fock spaces.
\begin{proposition}\label{propD01}
Let $u,\psi \in \mathcal{H}(\CC)$  such that $\psi$ is  not a constant,  $1\leq p, q\leq  \infty$,  and   $n\in \NN_0$. Then a bounded  $D_{(u,\psi,n )}: \mathcal{F}_p\to \mathcal{F}_q$   has closed range only if $p=q$.
\end{proposition}
The proof of the proposition will be given later  in Subsection~\ref{Proofp01}.

  We now   recall the notion of sampling sets for Banach spaces. The notion was introduced by
Ghatage, Zheng and Zorboska \cite{GZZ} as a tool  to study bounded below composition operators on  Bloch spaces. Since then,  it has been used to study both the  bounded below
and closed range  properties of  several    operators on spaces of analytic functions. There have been also various ways of  defining   the notion; see for example  \cite{GZZ,MPN}.  On  Fock spaces, we provide the following unified and general definition. Let  $1\leq p\leq  \infty,  k\in \NN_0$  and  $M$ be a non-empty subset of $\mathcal{F}_p$.   A  subset $S$ of $\CC$ is a $(p,k)$ sampling set  for $M$  if  there exists   a positive constant $\delta_k$ such that
\begin{align}
\label{sam}
\delta_k\|f\|_p \leq
\begin{cases}\sup_{z\in S} \frac{|f^{(k)}(z)|}{(1+|z|)^k}e^{-\frac{|z|^2}{2}}, \ \ p= \infty\\
\Big(\int_{S}\frac{|f^{(k)}(z)|^p}{(1+|z|)^{kp}} e^{-\frac{p}{2}|z|^2} dA(z)\Big)^{\frac{1}{p}}, \ p<\infty.
\end{cases}
\end{align} for all $f$ in $M$.
For each positive $\epsilon_n, n\in \NN_0$, we  also  define the sets
    \begin{align*}
   \Omega_{(u,\psi,n)}^{\epsilon_n}:= \{ z\in \CC : L_{(u,\psi,n)}(z) >\epsilon_n \}, \ \ \ G_{(u,\psi,n)}^{\epsilon_n}:= \psi\big( \Omega_{(u,\psi, n)}^{\epsilon_n}\big).
    \end{align*}
      We now state  our  next main result.
\begin{theorem} \label{closeddiff2}
  Let   $u, \psi\in \mathcal{H}(\CC)$ such that $\psi$ is not a constant,  $1\leq p\leq  \infty$,  $n\in \NN_0,$ and   $D_{(u,\psi,n)}$  be   bounded on $ \mathcal{F}_p $. Then  $D_{(u,\psi,n)}$ has closed range if and only if
   there exists $\epsilon_n>0$ such that  $ G_{(u,\psi,n)}^{\epsilon_n}$ is a $(p,n)$ sampling set for $ \mathcal{F}_{(p,n)}^0 $, where
  \begin{align*}
\mathcal{F}_{(p,n)}^0:= \{ f\in \mathcal{F}_p: f(0)= f'(0)= ...=f^{(n-1)}(0)=0 \}.
  \end{align*}
  \end{theorem}
\subsection{Surjective $D_{(u,\psi, n)}$  }
In this section we consider the question of  when  the operator  $D_{(u,\psi, n)}$ is surjective  on Fock spaces.  To state our result, we may first recall the notation of essential boundedness.   We  say a  non-zero  function $g$ in $\mathcal{H}(\CC)$  is  essentially bounded away from zero  if there exists a constant  $\delta >0$ such that   the measure of
    the set $\{z\in \CC: |g(z)|<\delta\}$ is zero.
   \begin{theorem} \label{thm88}
  Let   $u, \psi\in \mathcal{H}(\CC)$ such that $\psi$ is not a constant,  $1\leq p\leq \infty$, and $n\in \NN_0$. Let  $D_{(u,\psi, n)}$ be bounded on $\mathcal{F}_p $  and hence $\psi(z)= az+b$  for some $a, b\in \CC$ and $|a|\leq1$.
  Then the following statements are equivalent.
  \begin{enumerate}
  \item $  L_{(u,\psi,n)}$ is essentially bounded away from zero on $\CC$;
  \item    $D_{(u,\psi, n)}$   is surjective;
   \item  $|a|=1$.
      \end{enumerate}
\end{theorem}
We  close this section with   a word on notation.  The notion
 $U(z)\lesssim V(z)$ (or
equivalently $V(z)\gtrsim U(z)$) means that there is a constant
$C>0$ such that $U(z)\leq CV(z)$ holds for all $z$ in the set of
question. We write $U(z)\simeq V(z)$ if both $U(z)\lesssim V(z)$
and $V(z)\lesssim U(z)$.
\section{Proof of the  results} \label{mainproof}
In this section we present the proofs of the  results.  We begin by reminding the connection between the closed range problem  and the bounded below property of a linear operator on Banach spaces. Let $\mathcal{H}_1$ and $\mathcal{H}_2$ be two Banach spaces.  An operator  $T:\mathcal{H}_1 \to \mathcal{H}_2$ is said to be bounded below  if there exists a positive constant  $C$ such that $\|Tf\|_{\mathcal{H}_2}\geq C\|f\|_{\mathcal{H}_1}$ for every $f\in \mathcal{H}_1$.  As known from an application of the Open Mapping Theorem, an injective bounded operator on Banach spaces has  closed range if and only if it is  bounded below; see for example \cite[Theorem 2.5]{new}.\\
The operator  $D$ maps all constants  to the zero function and fails to be injective unless  its action is restricted   to  the spaces  modulo the constants or
\begin{align*}
\mathcal{F}_ {(m,p)}^0:= \big\{ f\in \mathcal{F}_ {(m,p)} : f(0)=0\big\}.\end{align*} In the latter  case, $D$ has closed range if and only if it is bounded below. On the other hand,
for each $f\in \mathcal{F}_ {(m,p)}$, the function $f- f(0)$ belongs to  $\mathcal{F}_ {(m,p)}^0$, and
$
D (f)= D(f-f(0))= f'. $
Thus,  $D$ has closed range on $\mathcal{F}_{(m,p)}^0$ if and only if it  has closed range on $\mathcal{F}_ {(m,p)}$.  For the sake of further referencing, we  record this  useful  observation below.
\begin{lemma}\label{lem3}
Let   $1\leq p,q\leq\infty $ and $m>0$.  Then  $D:\mathcal{F}_ {(m,p)}^0 \to \mathcal{F}_ {(m,q)} $ is bounded below if and only if    $D:\mathcal{F}_ {(m,p)}\to \mathcal{F}_ {(m,q)} $  has  closed range.
\end{lemma}
The lemma will be used in our next proof.
\subsection{Proof of Theorem~\ref{closeddiff1}}
Note that (iii) obviously implies (i).  Let us see  (i) $\Rightarrow$ (ii) and (ii) $\Rightarrow$ (iii).
 From \cite{Olivia2,TM3}, for  each $f\in \mathcal{F}_{(m,p)}$ we have
\begin{align}
\label{paley}
\|f\|_{(m, p)} \simeq \begin{cases}\Big(|f(0)|^p + \int_{\CC} \frac{|f'(z)|^p e^{-p|z|^m}}{ (1+ |z|)^{p(m-1)}} dA(z)\Big)^{\frac{1}{p}}, \ p<\infty\\
  |f(0)|+ \sup_{z\in \CC} \frac{|f'(z)| e^{-|z|^m}}{ (1+ |z|)^{(m-1)}}, \ \ p= \infty.
  \end{cases}
\end{align}
 We first consider the case when   either $p>q$  or  $p\leq q$ and $ m<2-\frac{pq}{pq+q-p}$.  Then  by Theorem~\ref{thmDD},  the operator is compact. It is known that a compact operator  has  closed range if and only if its range is finite dimensional. On the other hand, $D$ is injective on an infinite dimensional set  $\mathcal{F}_{(m,p)}^0$. Therefore, the operator has no closed range  in this case.

   Next, we consider  $p\leq q<\infty$ and $m=2-\frac{pq}{pq+q-p}$. For this case,  we  prove   the operator is not bounded below  unless $ p=q$. The norms of the monomials are estimated by   \begin{align}
\label{estmono}
\|z^n\|_{(p,  m)}
\simeq \bigg(\frac{n}{me }\bigg)^{\frac{n}{m}+\frac{2}{mp}-\frac{1}{2p}}.
\end{align} See \cite{Harutyunyan} for the details. For $p= \infty$, the corresponding estimate   becomes
\begin{align}
\label{estinf}
\|z^n\|_{(\infty, m)}\simeq \bigg(\frac{n}{me }\bigg)^{\frac{n}{m}}.
\end{align}
We may now use Lemma~\ref{lem3} and  suppose $p\leq q<\infty$ and the operator is bounded below. Then there exists a positive constant $\epsilon$ such that for all $n\in \NN$
\begin{align} \label{esti}
\|Dz^n\|_{(q,  m)}= n  \|z^{n-1}\|_{(q,  m)} \geq \epsilon \|z^{n}\|_{(p,  m)}.
\end{align} This  and  \eqref{estmono} imply
\begin{align}
\label{final}
\frac{n  \|z^{n-1}\|_{(q,  m)}}{\|z^{n}\|_{(p,  m)}} \simeq  \frac{n (n-1)^{\frac{n-1}{m}+\frac{2}{mq}-\frac{1}{2q}}}{n^{\frac{n}{m}+\frac{2}{mp}-\frac{1}{2p}}}\simeq
\frac{n^{\frac{m-1}{m}+\frac{2}{mq}-\frac{1}{2q}}}{n^{\frac{2}{mp}-\frac{1}{2p}}} \gtrsim \epsilon
\end{align} for all $n\in \NN$. Now setting $m=2-\frac{pq}{pq+q-p}$ and simplifying further, the relation in    \eqref{final} holds only when
$
n^{\frac{p-q}{2pq}} \gtrsim  \epsilon,
$ which implies  $p=q$  and hence  $m=2-\frac{pq}{pq+q-p}=1$.\\
Similarly,  if  $p\leq q=\infty$, then
\begin{align*}
\frac{n  \|z^{n-1}\|_{(\infty,  m)}}{\|z^{n}\|_{(p,  m)}} \simeq  \frac{n (n-1)^{\frac{n-1}{m}}}{n^{\frac{n}{m}+\frac{2}{mp}-\frac{1}{2p}}}
\simeq n^{1-\frac{1}{m}+\frac{1}{2p}-\frac{2}{mp}}=n^{-\frac{1}{2p}} \gtrsim \epsilon,
\end{align*} which implies  $p= q$ and hence $m=1$.

Next,  we show (ii) implies (iii). Let now  $p=q$ and  $m=1$. We need to show the range of the operator is $\mathcal{F}_{(1,p)}$. For each  $ f\in  \mathcal{F}_{(1,p)},$  consider the entire function
\begin{align*}h_f(z)= \int_{0}^z f(w)dw.\end{align*}  Applying  \eqref{paley},
\begin{align*}
\|h_f\|_{(1, p)}^p \simeq  \int_{\CC} |f(z)|^p e^{-p|z|} dA(z)\simeq\|f\|_{(1, p)}^p <\infty,
\end{align*} and hence $h_f\in \mathcal{F}_{(1,p)}$. Furthermore,   $Dh_f=  f$ and  completes the proof.
\subsection*{Remark 1}\label{remark} The same argument used above to show that a compact $D$  cannot have closed range on the Fock-type spaces  will be used  in the sequel for the operator $D_{(u,\psi, n)}$.
 \subsection{Proof of Theorem~\ref{thmorder}}
 First  note that for  $n\in \NN_0$ and $|z|\leq 1$,  using the  Cauchy integral formula and  \eqref{localest2} we have
\begin{align}
\label{est1}
| f^{(n)}(z)| \leq \frac{n!}{2\pi} \int_{|w-z|=1} \frac{|f(w)|}{|w-z|^{n+1}}|dw| \leq  n!\|f\|_p \max_{|w-z|=1} e^{\frac{|w|^2}{2}}\nonumber \\
\leq    n!e^{3/2} e^{\frac{|z|^2}{2}}\|f\|_p
\end{align}   for all $f\in \mathcal{F}_p$.  Similarly, for $|z|>1$,
\begin{align}
\label{est2}
| f^{(n)}(z)| \leq \frac{n!}{2\pi} \int_{|w-z|=1/|z|} \frac{|f(w)|}{|w-z|^{n+1}}|dw| \leq  n! |z|^n\|f\|_p \max_{|w-z|=1/|z|} e^{\frac{|w|^2}{2}}\nonumber \\
\leq n!e^{3/2}  |z|^n e^{\frac{|z|^2}{2}}\|f\|_p.
\end{align}
  Combining  \eqref{est1} and \eqref{est2}, we get
 \begin{align*}
 | f^{(n)}(z)| \leq n!e^{3/2} (1+|z|)^n e^{\frac{|z|^2}{2}}\|f\|_p,
 \end{align*}  and  hence
 \begin{align}
 \label{dd}
 |D_{(u, \psi, n)} f(z)|=  | u(z) f^{(n)}(\psi(z))| \leq  n!e^{3/2} |u(z)|(1+|\psi(z)|)^n e^{\frac{|\psi(z)|^2}{2}}\|f\|_p.
 \end{align}
  Suppose now that  $D_{(u, \psi, n)}$ is order bounded. Then there exists a positive  function $h_n\in L^q\big(\CC, e^{-\frac{q}{2}|z|^2}dA(z)\big)$  such that
  \begin{align*}
  | D_{(u, \psi, n)}f(z)|\leq h_n(z)
  \end{align*} for almost all $z\in \CC$ and $\|f\|_p\leq 1$.  Applying this inequality to the  normalized kernel functions $k_w, w\in \CC$, we get
 \begin{align*}
    | D_{(u, \psi, n)}k_w(z)|=  |u(z) k_w^{(n)}(\psi(z))|=   \big|u(z)\overline{w}^n \big| \big|e^{\overline{w}\psi(z)-|w|^2/2}\big| \leq h_n(z)
 \end{align*}   for almost all $z\in \CC$. For $w= \psi(z)$ in particular,
 \begin{align*}
\big|u(z)\overline{w}^n e^{\overline{w}\psi(z)-\frac{|w|^2}{2}}\big|= |u(z)\psi(z)^n | e^{\frac{|\psi(z)|^2}{2}}=L_{(u,\psi, n)}(z)e^{\frac{|z|^2}{2}}
\leq h_n(z)
 \end{align*} from which the necessity of the condition follows.\\
 To prove the converse,  setting \begin{align*}
  h_n(z):=  n! e^2 |u(z)| (1+|\psi(z)|)^n e^{\frac{|\psi(z)|^2}{2}}\simeq n! e^2  L_{(u,\psi, n)}(z)e^{\frac{|z|^2}{2}},  \end{align*} we observe that  the assumption on $L_{(u,\psi, n)}$
  implies  $h_n\in L^q\big(\CC, e^{-\frac{q}{2}|z|^2}dA(z)\big)$.  Furthermore, by  \eqref{dd}
\begin{align*}
 |D_{(u, \psi, n)} f(z)| \leq    n!e^{3/2}|u(z)| (1+|\psi(z)|)^{n} e^{\frac{|\psi(z)|^2}{2}} \leq h_n(z)
 \end{align*} for any  $ f\in\mathcal{F}_p$ such that $\|f\|_p \leq 1$, and completes the proof of the theorem.

 The next basic lemma  connects the closed range problem and boundedness from below of the operator $D_{(u,\psi, n)}$.
\begin{lemma}\label{lemcomp1}
Let $u,\psi\in \mathcal{H}(\CC)$ such that $\psi$ is not a constant,  $1\leq p\leq\infty $,  and $n\in \NN_0$.  Then a bounded
   $D_{(u,\psi, n)}: \mathcal{F}_p\to \mathcal{F}_p $  has  closed range if and only if the restriction operator  $D_{(u,\psi, n)}: \mathcal{F}_{(p,n)}^0\to \mathcal{F}_p $ is bounded from below, where \begin{align*}\mathcal{F}_{(p,n)}^0:= \big\{ f\in \mathcal{F}_p: f(0)= f'(0)= ...=f^{(n-1)}(0)=0 \big\}.\end{align*}
   \end{lemma}
   \emph{Proof}. Note  that  $D_{(u,\psi, n)}$ is  injective on $\mathcal{F}_{(p,n)}^0$ but not on   $\mathcal{F}_p$.  Thus,  as explained in the preceding section, $D_{(u,\psi, n)}$ has closed range on $\mathcal{F}_{(p,n)}^0$ if and only if   it is bounded below. On the other hand,   for each $f\in \mathcal{F}_p$, the function $f- S_n f$ belongs to $\mathcal{F}_{(p,n)}^0$,   where $S_n f$ refers to the first $n$ terms of  the  Taylor series expansion of the function $f$. Now,
    $
      D_{(u,\psi, n)} f= D_{(u,\psi, n)} (f-S_n f),
  $ from which  and the  connection between the closed range problem and boundedness below, the claim follows.
      \subsection{Proof of Proposition~\ref{propD01} }\label{Proofp01}
      Let  $1\leq p\leq \infty$.  From \cite{HU,UK}, the  estimate
 \begin{align}
  \label{HU}
  \|f\|_p \simeq
 \begin{cases}
\sum_{j=0}^{n-1}
|f^{(j)}(0)|
+ \Big(\int_{\CC} \frac{|f^{(n)}(z)|^p}{(1+|z|)^{np}}e^{-\frac{p}{2}|z|^2}dA(z)\Big)^{\frac{1}{p}}, \ p<\infty\\
 \sum_{j=0}^{n-1}
|f^{(j)}(0)|+
 \sup_{z\in \CC} \frac{ |f^{(n)}(z)|}{ (1+|z|)^{n}}e^{-\frac{1}{2}|z|^2}, \ p= \infty
 \end{cases}
\end{align} holds. We will appeal to  this estimate several times in the sequel. \\
 Let us consider  first the case $p<q<\infty$  and assume $D_{(u,\psi, n)}:\mathcal{F}_p\to \mathcal{F}_q$ has closed range. By Lemma~\ref{lemcomp1},  the operator is bounded  below
 on $\mathcal{F}_{(p,n)}^0$. We consider the  sequence of the monomials  $f_k(z)= z^k, k\in \NN_0, k\geq n$. Using  Stirling's  approximation formula again
\begin{align}
\label{normp}
 \| f_k\|_p^p= p \int_{0}^\infty r^{kp+1} e^{-pr^2/2} dr
= \big(1/p\big)^{kp/2} \Gamma\big((kp+2)/2\big) \simeq \big(k/e\big)^{\frac{kp}{2}} \sqrt{k}.
\end{align} See also \cite[p. 40]{Zhu}.
  Now  applying  \eqref{HU} and Theorem~\ref{thm3},
\begin{align*}
\|D_{(u,\psi,n )} f_k\|_q^q =  \frac{q}{2\pi} \int_{\CC} |u(z)|^q  |f_k^{(n)}(az+b)|^{q} e^{-\frac{q}{2}|z|^2} dA(z) \quad \quad\quad \quad \quad \quad\nonumber\\
=\frac{q}{2\pi} \int_{\CC} |u(z)|^q (1+|az+b|)^{nq} \frac{|f_k^{(n)}(az+b)|^{q}}{(1+|az+b|)^{nq}} e^{-\frac{q}{2}|z|^2} dA(z) \quad  \nonumber\\
\simeq \int_{\CC}  L_{(u,\psi,n)}^q(z) \frac{|f_k^{(n)}(az+b)|^{q}}{(1+|az+b|)^{nq}} e^{-\frac{q}{2}|az+b|^2} dA(z) \quad \quad \\
\lesssim L^q_n \int_{\CC} \frac{|f_k^{(n)}(az+b)|^{q}}{(1+|az+b|)^{nq}} e^{-\frac{q}{2}|az+b|^2} dA(z)
\lesssim  L^q_n \|f_k\|_q^q \lesssim  \|f_k\|_q^q,  \end{align*}
where $L_n$ is a constant as in  Theorem~\ref{thm3}.
This and  boundedness below imply there exists $\epsilon_n >0$ for which
\begin{align}
\label{equal}
 \|f_k\|_q \geq \epsilon_n  \|f_k\|_p.
\end{align}
  Now, applying  \eqref{normp}, the estimate in \eqref{equal} holds only if
$k^{\frac{1}{2q}-\frac{1}{2p}} \geq \epsilon_n$ for all $k\in \NN_0, k\geq n$. This  gives a   contradiction when $k\to \infty$.\\
Similarly, for $p<q= \infty$, we have
$
 \| f_k\|_\infty=\big(k/e\big)^{k/2}.
$
 By  \eqref{HU}  and Theorem~\ref{thm3},
\begin{align*}
\|D_{(u,\psi,n)} f_k\|_\infty \simeq \sup_{z\in \CC} L_{(u,\psi,n)}(z)\frac{|f_k^{(n)}(az+b)|}{(1+|az+b|)^{n}} e^{-\frac{1}{2}|az+b|^2} \lesssim L_n\|f_k\|_\infty.
\end{align*} Therefore,
\begin{align*}
 \Big(\frac{k}{e}\Big)^{\frac{k}{2}}=\|f_k\|_\infty \geq \epsilon_n  \|f_k\|_p \simeq \epsilon_n \Big(\frac{k}{e}\Big)^{\frac{k}{2}} k^{\frac{1}{2p}}
 \end{align*} for some $\epsilon_n >0.$ This gives a contradiction when $k\to \infty $ again.\\
 If $p>q$, then by Theorem~\ref{thm3}, $D_{(u,\psi,n )}:\mathcal{F}_p\to \mathcal{F}_q$ is compact and  injective  on  an infinite dimensional space $\mathcal{F}_{(p,n)}^0$. By Remark 1, it follows that  the range of the operator  cannot be  closed.
 \subsection{Proof of Theorem~\ref{closeddiff2}} \label{sub1}
     By Lemma~\ref{lemcomp1}, it is enough to show that    $D_{(u,\psi,n)}$   has closed range on $\mathcal{F}_{(p,n)}^0$
 if and only if  there exists a constant   $\epsilon_n>0$ such that  $G_{(u,\psi, n)}^{\epsilon_n}$ is a $(p,n)$ sampling set for $\mathcal{F}_{(p,n)}^0$.

Suppose $G_{(u,\psi, n)}^{\epsilon_n}$ is a $(p,n)$  sampling set for some $\epsilon_n>0$. If $  p<\infty$, then  there exists  $\delta_n>0$ such that for each $f \in \mathcal{F}_{(p,n)}^0$,
\begin{align}
\label{bb1}
\delta_n\|f\|_p^p \leq  \int_{G_{(u,\psi, n)}^{\epsilon_n}} \frac{ |f^{(n)}(z)|^p}{(1+|z|)^{np}} e^{-\frac{p}{2}|z|^2} dA(z).
\end{align} It follows that
\begin{align*}
\|D_{(u,\psi, n)}f\|_p^p= \frac{p}{2\pi}\int_{\CC} \big| u(z) f^{(n)}(\psi(z))\big|^p e^{-\frac{p}{2}|z|^2} dA(z)\quad \quad \quad \quad \quad \quad \quad \quad \quad \\
 \geq  \frac{p}{2\pi}  \int_{\Omega_{(u,\psi, n)}^{\epsilon_{n}}}  \big| u(z) f^{(n)}(\psi(z))\big|^p e^{-\frac{p}{2}|z|^2}dA(z)\quad \quad \quad \quad \quad \quad \quad\\
\geq \frac{p}{2\pi} \int_{\Omega_{(u,\psi, n)}^{\epsilon_n}} L_{(u,\psi, n)}^p(z) \frac{| f^{(n)}(\psi(z))|^p}{(1+|\psi(z)|)^{np}} e^{-\frac{p}{2}|\psi(z)|^2} dA(z).\quad \quad \quad
\end{align*}
By \eqref{bb1}, the last right-hand  integral above is bounded below by
\begin{align}
\label{below1}
 \frac{p \epsilon_n^p}{2\pi |a|^2}  \int_{G_{(u,\psi, n)}^{\epsilon_n}}  \frac{| f^{(n)}(z)|^p}{(1+|z|)^{np}} e^{-\frac{p}{2}|z|^2} dA(z)
\geq  \frac{p \epsilon_n^p}{2\pi |a|^2}   \delta_n   \|f\|_p^p.\quad \quad \quad \quad
\end{align}
Similarly, for   $p=\infty$,  there exists $\delta_n$ such that  for each $f \in \mathcal{F}_{(p,n)}^0$,
\begin{align}
\label{below200}
\quad \|D_{(u,\psi, n)}f\|_\infty = \sup_{z\in\CC} \big| u(z) f^{(n)}(\psi(z))\big| e^{-\frac{1}{2}|z|^2} \quad \quad \quad \quad \quad \quad \quad \quad \quad \quad \quad \nonumber\\
 \geq \sup_{z\in \Omega_{(u,\psi, n)}^{\epsilon_n}} L_{(u,\psi, n)}(z) (1+|\psi(z)|)^{-n} | f^{(n)}(\psi(z))| e^{-\frac{1}{2}|\psi(z)|^2}  \nonumber \\
\geq  \epsilon_n \sup_{z\in G_{(u,\psi, n)}^{\epsilon_n}} (1+|z|)^{-n}| f^{(n)}(z)| e^{-\frac{1}{2}|z|^2}
\gtrsim \delta_n\epsilon_n  \|f\|_\infty.\quad
\end{align}
From \eqref{below1}, \eqref{below200} and Lemma~\ref{lemcomp1}, the sufficiency of the condition follows.\\
Conversely,  let $n$ be fixed and suppose    $G_{(u,\psi,n )}^{\epsilon}$ is not a  $(p,n)$  sampling set for  each  $\epsilon>0$.   Let $p<\infty$. Then  there exists a unit norm sequence
$(f_{k})_{ k\in \NN}$ in $\mathcal{F}_{(p,n)}^0$ such that  \begin{align*}
\int_{G_{(u,\psi,n)}^{1/k}} \frac{|f_{k}^{(n)}(z)|^p}{(1+|z|)^{np}} e^{-\frac{p}{2}|z|^2} dA(z) \to 0 \   \text{as}  \  k\to \infty.
\end{align*}
It follows that
\begin{align*}
\|D_{(u,\psi,n)}f_k\|_p^p= \frac{p}{2\pi}\int_{\CC}  |u (z)|^p |f_k^{(n)}(\psi(z))|^p e^{-\frac{p}{2}|z|^2} dA(z)\quad \quad \quad \\
=  \frac{p}{2\pi}\int_{\CC} L_{(u,\psi, n)}^p(z)\frac{ |f_k^{(n)}(\psi(z))|^p}{(1+|\psi(z)|)^{np}}e^{-\frac{p}{2}|\psi(z)|^2} dA(z)
= I_{1n}+ I_{2n},
\end{align*}
where we set
\begin{align*}
 I_{1n}  =\frac{p}{2\pi}\int_{\Omega_{(u,\psi,n)}^{1/k}}L_{(u,\psi, n)}^p(z)\frac{ |f_k^{(n)}(\psi(z))|^p}{(1+|\psi(z)|)^{np}}e^{-\frac{p}{2}|\psi(z)|^2} dA(z)\\
  \leq \frac{pL_n^p}{2\pi}  \int_{\Omega_{(u,\psi,n)}^{1/k}}\frac{ |f_k^{(n)}(\psi(z))|^p}{(1+|\psi(z)|)^{np}}e^{-\frac{p}{2}|\psi(z)|^2} dA(z)\quad\quad \quad \quad \\
  =\frac{pL_n^p}{2\pi |a|^2} \int_{G_{(u,\psi,n)}^{1/k}}\frac{ |f_k^{(n)}(z)|^p}{(1+|z|)^{np}}e^{-\frac{p}{2}|z|^2} dA(z)\to 0 \quad\quad\quad
 \end{align*} as $ k\to \infty.$
To estimate the remaining integral, we eventually apply \eqref{HU} and
\begin{align*}
 I_{2n} \lesssim \int_{\CC\setminus \Omega_{(u,\psi,n)}^{1/k}}  L_{(u,\psi,n)}^p(z) \frac{|f_k^{(n)}(\psi(z))|^p}{(1+|\psi(z)|)^{np}} e^{-\frac{p}{2}|\psi(z)|^2} dA(z) \quad \\
\leq \frac{1}{k^p} \int_{\CC\setminus \Omega_{(u,\psi, n)}^{1/k}}   \frac{|f_k^{(n)}(\psi(z))|^p}{(1+|\psi(z)|)^{np}} e^{-\frac{p}{2}|\psi(z)|^2} dA(z) \quad  \quad \quad \quad\\
\leq  \frac{1}{k^p} \int_{\CC} \frac{|f_k^{(n)}(\psi(z))|^p}{(1+|\psi(z)|)^{np}} e^{-\frac{p}{2}|\psi(z)|^2} dA(z)\quad\quad \quad \quad \quad\quad \quad \\
\simeq \frac{\|f_k\|_p^p}{k^p}  \simeq \frac{1}{k^p} \to 0  \  \text{as}  \ k\to \infty. \quad\quad \quad \quad\quad\quad \quad \quad  \quad  \quad
\end{align*}  This   contradicts the assumption that the operator is bounded below.

Next, consider $p=\infty$ and suppose $D_{(u,\psi, n)}$ is bounded  below. Then there exists a constant  $\delta_n>0$ such that  for each $f \in \mathcal{F}_{(\infty,n)}^0$
\begin{align*}
\sup_{z\in \CC}|u(z)| |f^{(n)}(\psi(z))|e^{-\frac{|z|^2}{2}} \geq \delta_n \| f\|_\infty.
\end{align*}
Then, by definition of supremum for each $f$  there exists $w_f\in \CC$ such that
\begin{align}
\label{small}
|u(w_f)| |f^{(n)}(\psi(w_f))|e^{-\frac{|w_f|^2}{2}} > \frac{\delta_n}{2} \| f\|_\infty.
\end{align}
On the other hand, by \eqref{HU} again
\begin{align*}
|u(w_f)| |f^{(n)}(\psi(w_f))|e^{-\frac{|w_f|^2}{2}}\lesssim  L_{(u,\psi, n)} (w_f) \frac{|f^{(n)}(\psi(w_f))|}{(1+|\psi(w_f)|)^{n}}e^{-\frac{|\psi(w_f)|^2}{2}}\\
\leq L_{(u,\psi, n)} (w_f) \| f\|_\infty  \quad\quad \quad \quad \quad\quad  \quad\quad
\end{align*} and with \eqref{small}, we deduce
\begin{align*}
 L_{(u,\psi, n)} (w_f) > \frac{\delta_n}{2}.
\end{align*}
Setting $\epsilon_n= \delta_n/2$, we observe that $w_f\in \Omega_{(u,\psi, n)}^{\epsilon_n}$ and using \eqref{small}
\begin{align*}
\| f\|_\infty \leq \frac{2 |u(w_f)|}{\delta_n} |f^{(n)}(\psi(w_f))|e^{-\frac{|w_f|^2}{2}} \quad \quad \quad \quad \quad \quad \quad \quad \quad \quad\quad \\
\lesssim \frac{2L_{(u,\psi, n)}(w_f)}{\delta_n}  \frac{|f^{(n)}(\psi(w_f))|}{(1+|\psi(w_f)|)^{n}}e^{-\frac{|\psi(w_f)|^2}{2}}\\
 \lesssim \frac{2L_{n}}{\delta_n} \frac{|f^{(n)}(\psi(w_f))|}{(1+|\psi(w_f)|)^{n}}e^{-\frac{|\psi(w_f)|^2}{2}}\quad\quad \quad \quad  \\
 \leq \frac{2L_{n}}{\delta_n}  \sup_{z\in G_{(u,\psi,n)}^{\epsilon_n}} \frac{|f^{(n)}(z)|}{(1+|z|)^{n}}e^{-\frac{|z|^2}{2}}\quad
\end{align*} and completes the proof.
      \subsection{Proof of Theorem~\ref{thm88} }
Let $p<\infty$.  We prove first  the implication (i)$\Rightarrow $(ii) and suppose $\gamma_n$ is an essential lower bound for  $ L_{(u,\psi, n)}$.  Then for each  $ f\in \mathcal{F}_p$, consider the function
  \begin{align*}
  h_f(z)= \begin{cases}
   g(z)\ \  \ \ \ u(z)\neq 0\\
  \lim_{w \to z } g(w),\ \  \ \           u(z)=0,
    \end{cases}
  \end{align*}
  where we set
  \begin{align*}
  g(z)= \int_0^z \int_0^{z_1}  \int_0^{z_2}...  \int_0^{z_{n-1}} \frac{f(\psi^{-1}(w))}{u(\psi^{-1}(w))} dA(w) dA(z_{n-1})... dA(z_2) dA(z_1).
  \end{align*}
  Clearly, $D_{(u,\psi, n)}h_f=  f$.
    Since $ u $ is entire and vanishes  at most in a set of measure zero, we estimate
    \begin{align*}
  \| h_f\|_p^p \simeq  \int_{\CC} \frac{ |g^{(n)}(z)|^p}{(1+|z|)^{np}}  e^{-\frac{p}{2}|z|^2}dA(z) \simeq  \int_{\CC} \frac{|f(\psi^{-1}(z))|^p e^{-\frac{p}{2}|z|^2}}{(1+|z|)^{np}|u(\psi^{-1}(z))|^p}dA(z)\\
\leq \int_{\CC} \frac{|f(z)|^p e^{-\frac{p}{2}|\psi(z)|^2}}{(1+|\psi(z)|)^{np}|u(z)|^p}dA(z)\quad \\
\lesssim \int_{\CC} L_{(u,\psi, n)}^{-p}(z)  |f(z)|^p e^{-\frac{p}{2}|z|^2}dA(z) \\
\leq 2\pi\gamma_n^{-p}p^{-1} \|f\|_p^p <\infty, \quad \quad \quad\quad \quad
    \end{align*} from which the statement in (ii) follows.

 Next, we prove $(ii) \Rightarrow (iii)$, and suppose on the contrary $|a|<1$.  The surjectivity of the operator implies there exists some   $ f\in \mathcal{F}_p$ such that
    $1=  D_{(u,\psi,n)}f$. It follows that $u$ has no zeros in  $\CC$, and $u\psi^n$ has a zero set of measure zero which  does not affect our integral approach below in \eqref{surj}. Thus, we can assume that
    $u\psi^n$  is non-vanishing.     Then by  Lemma~\ref{lemnew},   it follows that
  $  u\psi^n(z)= e^{a_{0n}+a_{1n}z+ a_{2n}z^2},
   $ for some constants $a_{0n}, a_{1n},  a_{2n} \in \CC $  such that $|a_{2n}|\leq \frac{1-|a|^2}{2}$.  Using this, we write
   \begin{align*}
L_{(u, \psi, n)} (z)= |u(z)\psi^n(z)|e^{\frac{1}{2}(|az+b|^2-|z|^2)}\quad \quad \quad \quad \quad \quad \quad \quad \quad \quad \nonumber \\
=  e^{\Re(a_{0n}+a_{1n}z+ a_{2n}z^2)}e^{\frac{1}{2}(|az+b|^2-|z|^2)} \quad \quad \quad \quad \quad \quad \quad \quad  \\
= C e^{\Re((a_{1n}+a\overline{b})z)+\Re(a_{2n}z^2)+\frac{|a|^2-1}{2}|z|^2} \quad \quad \quad\quad \quad \quad \quad \quad
\end{align*} for all $z\in \CC ,$  where   $C= e^{\Re(a_{0n}) + \frac{|b|^2}{2}}.$
 By surjectivity, for each $ h\in \mathcal{F}_p$, there exists some $ f\in \mathcal{F}_p$  such that $D_{(u,\psi,n)}f(z)= u(z) f^{(n)}( az+b)= h(z)$ for all $z\in \CC$.  This implies
 \begin{align}
 \label{surj}
 \int_{\CC} \frac{|f^{(n)}( az+b)|^p   e^{-\frac{p}{2}|az+b|^2}}{(1+|\psi(z)|)^{np}} dA(z)=  \int_{\CC} \frac{|h(z)|^p e^{-\frac{p}{2}|z|^2}}{L^{p}_{(u, \psi, n)} (z)} dA(z).
 \end{align}
 By \eqref{HU}, the right-hand integral in \eqref{surj} should be finite for each $ h\in \mathcal{F}_p$. The plan is now  to show the existence of some functions   $h $ in the space for which this integral diverges.

Now, if $|a_{2n}|< \frac{1-|a|^2}{2}$, then the operator is compact and by  Remark~1, its range is not closed and the operator is not surjective. Thus, we set  $|a_{2n}|= \frac{1-|a|^2}{2}$ and consider the following two cases following Lemma~\ref{lemnew}.\\
\textbf{Case 1.}  For $a_{1n}+a\overline{b}=0$,  we have
\begin{align}
\label{check}
L_{(u, \psi, n)} (z)=  C e^{\Re(a_{2n}z^2)+\frac{|a|^2-1}{2}|z|^2}.
\end{align}
We may  also  write $a_{2n} = |a_{2n}| e^{-2i\theta_{2n}} $, where $0\leq \theta_{2n}<\pi$.\\
 Replacing
 $z $ by $e^{i\theta_{2n}}w$   in \eqref{check}
 \begin{align}
 \label{lett}
 L_{(u, \psi, n)} (e^{i\theta_{2n}}w)=  C  e^{ \frac{1-|a|^2}{2}(\Re(w^2)-|w|^2)}
\end{align} for all $w\in \CC$.  Setting $w= x+iy$,  the relation in \eqref{lett} implies
\begin{align}
\label{lett0}
\frac{ e^{-\frac{1}{2}|e^{i\theta_{2n}}w|^2}}{L_{(u, \psi, n)} (e^{i\theta_{2n}}w)} = \frac{1}{C} e^{(\frac{1}{2}-|a|^2)y^2- \frac{1}{2}x^2}.
\end{align}
from which we observe that the integral in \eqref{surj} diverges for every nonzero constant function $h$  in the space whenever $|a|< \frac{1}{\sqrt{2}}$. Thus,  the question is  when $|a|\geq  \frac{1}{\sqrt{2}}$.  We may  consider
a function $h(z)= h_\alpha(z)= e^{\alpha z^2}$  where $\alpha$ is a real number  and  $|\alpha| < \frac{1}{2}$. A suitable $\alpha$ will be chosen later.
A straightforward calculation using \eqref{lett0} gives
\begin{align}
\label{int}
\int_{-\infty}^\infty \int_{-\infty}^\infty
\frac{|h_\alpha(e^{i\theta_{2n}}(x+iy))|^p e^{-\frac{p}{2}|e^{i\theta_{2n}}(x+iy)|^2}}{L^{p}_{(u, \psi, n)} (e^{i\theta_{2n}}(x+iy))} dx dy=\frac{1}{C}\int_{-\infty}^\infty \int_{-\infty}^\infty  e^{p\beta(x,y)} dx dy,
\end{align} where
\begin{align*}
\beta(x,y):= \big(\alpha \cos(2\theta_{2n})-\frac{1}{2}\big) x^2 +  \big(\frac{1}{2}-|a|^2-\alpha \cos(2\theta_{2n})\big) y^2\nonumber\\
-2\alpha xy \sin(2\theta_{2n}).
\end{align*}
 Since $|\alpha| <1/2$,  it holds that $ \alpha \cos(2\theta_{2n})-\frac{1}{2}< 0$. Integrating  with  respect to $x$,
\begin{align*}
\int_{-\infty} ^\infty  e^{ p(\alpha \cos(2\theta_{2n})-\frac{1}{2}) x^2 - 2p\alpha y \sin(2\theta_{2n})x} dx
= \int_{-\infty} ^\infty  e^{ -p(\frac{1}{2}-\alpha \cos(2\theta_{2n})) x^2 - 2p\alpha y \sin(2\theta_{2n})x} dx \quad \nonumber\\
=\sqrt{\pi}\Big(\frac{p}{2}-p\alpha \cos(2\theta_{2n})\Big)^{-\frac{1}{2}} e^{\frac{\alpha^2 y^2p^2 \sin^2(2\theta_{2n})}{\frac{p}{2}-p\alpha \cos(2\theta_{2n})}}.\quad \quad
\end{align*}
Taking this in \eqref{int}, the coefficient of  $y^2$ becomes
\begin{align}
\label{choice}
\frac{p}{2}-p|a|^2-p\alpha \cos(2\theta_{2n})+\frac{\alpha^2 p^2 \sin^2(2\theta_{2n})}{\frac{p}{2}-p\alpha \cos(2\theta_{2n})}\quad \quad \quad \quad \quad \quad\quad \quad \quad\nonumber\\
=\frac{p\alpha^2 +\frac{p}{4}-\frac{p|a|^2}{2}+p\alpha \cos(2\theta_{2n})(|a|^2-1)}{\frac{1}{2}-\alpha \cos(2\theta_{2n})}.
\end{align}
Now, if  $\cos(2\theta_{2n}) \leq 0$,  we choose a positive  $\alpha $ such that
\begin{align}
\label{neg}
\frac{1}{4} >\alpha^2 >\frac{|a|^2}{2}- \frac{1}{4}.
\end{align}
 Note that such a choice is possible since $|a|<1$.  For such $\alpha$, the expression in
 \eqref{choice} is nonnegative and hence the integral in \eqref{int} diverges.  On the other hand, if  $\cos(2\theta_{2n}) > 0$,  we can  choose a  negative  $\alpha $ such that  \eqref{neg} holds and hence the integral in \eqref{int} diverges again.

\textbf{Case 2}.  Let   $a_{1n}+a\overline{b}\neq 0$ and \begin{align*}a_{2n}= - \frac{(1-|a|^2)(a_{1n}+a\overline{b})^2}{2|a_{1n}+a\overline{b}|^2}.\end{align*}
Using this and  $a_{2n} = |a_{2n}| e^{-2i\theta_{2n}} $ as above, we obtain
 \begin{align*}(a_{1n}+a\overline{b} )e^{i\theta_{2n}} =\pm i |a_{1n}+a\overline{b}|,
 \end{align*}  which is a purely imaginary number. Setting  $(a_{1n}+a\overline{b} )e^{i\theta_{2n}}= iy_n$ for
some $y_n\in \RR$,  $w= x+iy$, and $ z=e^{i\theta_{2n}}w$
\begin{align*}
L_{(u, \psi, n)} (e^{i\theta_{2n}}w)=  C  e^{ -y_ny- (1-|a|^2)y^2}
\end{align*} and hence
\begin{align*}
\frac{ e^{-\frac{1}{2}|e^{i\theta_{2n}}w|^2}}{L_{(u, \psi, n)} (e^{i\theta_{2n}}w)} = \frac{1}{C} e^{y_ny+(\frac{1}{2}-|a|^2)y^2- \frac{1}{2}x^2}.
\end{align*}
This shows that if  $|a|\leq \frac{1}{\sqrt{2}}$, then  the integral in \eqref{surj} diverges for every nonzero constant function $h$ in the space again. For the rest,
 we consider  a function $ h_\alpha(z)= e^{\alpha z^2}$  and argue  exactly in the same  way as above.

It remains to show  (iii)$\Rightarrow$ (i).  For  $|a|=1$, by \eqref{kernform} we have
    \begin{align*}
  L_{(u, \psi, n)} (z)=  |u(z)\psi^n(z)| e^{\frac{1}{2}(|az+b|^2-|z|^2)}= |b^nu(0)| e^{\frac{|b|^2}{2}}.
  \end{align*}
  Note that if $b^nu(0)=0$,  then the function $u$ vanishes since $\psi$ is entire and nonzero.  Hence  $b^nu(0)\neq0$ and $L_{(u, \psi, n)}$ is bounded away from zero.

  For $p= \infty$, we simply replace the integral argument in the proof  by the  supremum.
\subsection*{Acknowledgement}
I would like to thank the referee for constructive comments and suggestions which helped improve  the  manuscript.


\begin{thebibliography}{BRSHZE}
\bibitem{new} Y. Abramovich and C. Aliprantis,  An invitation to operator theory.  Graduate studies in Mathematics, vol. 50. American Mathematical Society, Providence (2002).

\bibitem{Maria} M.  Beltr\'{a}n, Dynamics of differentiation and integration operators on weighted spaces of entire functions, Studia Matematica,  \textbf{221} (2014), 35--60.

\bibitem{Bon} J. Bonet, A. Bonilla, Chaos of the differentiation operator on weighted
Banach spaces of entire functions, Complex Anal. Oper. Theory,  \textbf{7 }(2013),no. 1, 33--42.

\bibitem{BMM} J. Bonet, T. Mengestie, and  M. Worku,  Dynamics of the Volterra-type integral and differentiation operators on generalized Fock spaces,  Results Math 74, \textbf{197} (2019). https://doi.org/10.1007/s00025-019-1123-7

\bibitem{CG} T. Carroll and  C. Gilmore,  Weighted composition operators on Fock Spaces and their dynamics, J. Math. Anal. Appl.,
\textbf{502} (2021),  125234.

\bibitem{Olivia2} O. Constantin and A.  Persson, The spectrum of Volterra type integration operators on generalized  Fock spaces,
 Bull. London Math. Soc.,   \textbf{47} (2015),  958--963.

\bibitem{JDHJ} J. Diestel and H. Jarchow, Absolutely Summing Operators. Cambridge University Press, Cambridge (1995).

\bibitem{TD} T. Domenig,  Order bounded and p-summing composition operators. Contemp. Math.,  \textbf{213} (1998), 27--41.

\bibitem{GZZ} P. Ghatage, D. Zheng, and N. Zorboska, Sampling sets and closed range composition operators on
the Bloch space. Proc. Amer.Math. Soc.,  \textbf{133} (2005), 1371--1377.


 \bibitem{Harutyunyan}  A. Harutyunyan and W. Lusky,
On the boundedness of the differentiation operator
between weighted spaces of holomorphic functions, Studia Math.,  \textbf{184} (2008), 233--247.

\bibitem{HP} R.  Hibschweiler  and N. Portnoy,  Composition followed by differentiation between Bergman and Hardy spaces. Rocky Mt. J. Math., \textbf{35} (2005),  843--855.
\bibitem {HHJ} H. Hunziker and H.Jarchow, Composition operators which improve integrability. Math. Nachr., \textbf{152} (1991), 83--99.


\bibitem{HU} Z. J. Hu, Equivalent norms on Fock spaces with some application to extended Ces\`{a}ro
operators, Proc Amer Math Soc., \textbf{141} (2013),  2829--2840.

\bibitem{TMD} T.  Mengestie,   A note on the differentiation  operator on generalized Fock spaces,  J. Math. Anal. Appl., \textbf{458} (2018), 937--948.

\bibitem{TMDV} T. Mengestie, On the differential and Volterra-type integral operators on Fock-type spaces, Revista de la Unión Matemática Argentina, \textbf{63} (2022),   379--395.

\bibitem{TM5} T. Mengestie, Spectral properties  of Volterra-type integral operators on Fock--Sobolev  spaces,   J. Korean Math. Soc.,  \textbf{54} (2017),  1801--1816.

\bibitem{TM3} T. Mengestie and S.  Ueki,  Integral, differential  and multiplication operators on weighted Fock spaces, Complex Anal. Oper. Theory,  \textbf{13} (2019), 935--958.

\bibitem{MPN} M. M. Pirsastech, N. Eghbali, and A. H. Sanatpour,  Closed range properties of Li–Stević integral-type operators between Bloch-type spaces and their essential norms, Turk. J. Math.,
 \textbf{ 42} (2018),  3101-- 3116.

\bibitem{AKS} A. K. Sharma,  On order bounded weighted composition operators between Dirichlet spaces. Positivity, \textbf{21} (2017), 1213--1221.

\bibitem{UK} S. Ueki, Higher order derivative characterization for Fock-type spaces, Integr. Equ.
Oper. Theory, \textbf{84} (2016), 89--104.

\bibitem{EW} E. Wolf, Order bounded weighted composition operators, J. Aust. Math. Soc., \textbf{93} (2012), 333--343.

\bibitem{MW} M. Worku, On generalized weighted composition operators acting between Fock-type spaces, Quaestiones Mathematicae, \textbf{45} (2022), 1317--1331.

\bibitem{Zhu}  K. Zhu, Analysis on Fock spaces. Graduate Texts in Mathematics, 263. Springer, New York, 2012. x+344.
\end{thebibliography}
\end{document}